\documentclass[12pt,letterpaper]{amsart}
\usepackage{amsmath,amssymb,amsfonts}
\usepackage[all]{xy}
\usepackage{graphicx}
\usepackage{caption}
\usepackage{enumerate}
\usepackage{xcolor}

\usepackage{appendix}
\usepackage{a4wide}
\usepackage{hyperref}

\newtheorem*{thm}{Theorem}

\theoremstyle{definition}
\newtheorem*{definition}{Definition}

\newtheorem*{rem}{Remark}

    \DeclareFontFamily{U}{wncy}{}
    \DeclareFontShape{U}{wncy}{m}{n}{<->wncyr10}{}
    \DeclareSymbolFont{mcy}{U}{wncy}{m}{n}
    \DeclareMathSymbol{\Sha}{\mathord}{mcy}{"58}

\numberwithin{equation}{section}

\DeclareSymbolFont{bbold}{U}{bbold}{m}{n}
\DeclareSymbolFontAlphabet{\mathbbold}{bbold}

\DeclareSymbolFont{bbold}{U}{bbold}{m}{n}
\DeclareSymbolFontAlphabet{\mathbbold}{bbold}


\title{G-circulant Matrices and the Classical Maschke Theorem}
\author{Jon Merzel}
\address{Department of Mathematics, Soka University of America, 1 University Drive, Aliso Viejo, CA 92656}
\email{jmerzel@soka.edu}

\date{\today}

\begin{document}
\maketitle
\begin{abstract}
In this note, we use the isomorphism of the ring of $G$-circulant matrices over
a field $k$ with the group ring $k[G]$ to derive a very short proof of the Classical Maschke Theorem.
    
\end{abstract}

\bigskip 

\begin{center} 
Dedicated to J\'{a}n Min\'{a}\v{c} for his 70th birthday
\end{center}
\bigskip

By ``the Classical Maschke Theorem" we mean:

\begin{thm}
Let $k$ be a field, $G$ \ a finite group whose order $n$ is not divisible by
the characteristic of $k$. \ Then the group ring $k[G]$ is semisimple.
\end{thm}

Of course, as $k[G]$ is finite dimensional over $k$, one needs to show only
that the Jacobson radical $rad(k[G])=0$. \ We want to use the following basic
fact about the Jacobson radical (see e.g. \cite{[Lam]}, Theorem 4.12): \ If $R$ is a
left artinian ring, then $rad(R)$ is a nil (in fact nilpotent) ideal. (An ideal $I$ is nil if each element of $I$ is nilpotent; $I$ is nilpotent if $I^k={0}$ for some positive integer $k$.)  Of course, this is then
the case when $R$ is the group algebra $k[G]$, $k$ a field, $G$ a finite group.

What we actually need here is that each element of $rad(k[G])$ is nilpotent;
to keep things elementary, we  give an ad hoc proof of a (sufficiently
general) special case of the cited result. \ Let $R$ be a finite dimensional
(associative) algebra over the field $k$, and $\alpha \in rad(R)$. \ Then $%
u+\alpha b$ is a unit in $R$ for any $b\in R$ and any unit $u\in R$, so if $%
f(x)\in k[x]$ and $f(0)\neq 0$ then $f(\alpha )$ is a unit in $R$. \ Since $R
$ is finite dimensional, $g(\alpha )=0$ for some nonzero polynomial $g(x)\in
k[x]$; write $g(x)=x^{m}f(x)$ with $f(0)\neq 0$. \ As $f(\alpha )$ is a
unit, we must have $m>0$ and $\alpha ^{m}=0$, so $\alpha $ is nilpotent.\

We now give a brief explicit description of the right regular representation of $G$.

\begin{definition}
Let $G=\{g_{i}\}_{i=1}^{n}$ be a finite group written multiplicatively, and $
k$ a field. \ For each $g\in G$ we define the $n\times n$ matrix $A^{g}\in
M_{n,n}(k)$ by $A_{i,j}^{g}=\delta _{g_{i}g,g_{j}}$ (Kronecker $\delta $).
\end{definition}

Simple computations show that $A^{g}$ is a permutation matrix, that the
matrices $A^{g_{1}},\cdots ,A^{g_{n}}$ are linearly independent (their
nonzero entries occur in pairwise disjoint sets of indices), and that $%
A^{g}A^{h}=A^{gh}$. \ Thus $\{A^{g}|~g\in G\}$ under matrix multiplication
forms a group isomorphic to $G$ and $\{\sum\limits_{i}a_{i}A^{g_{i}}|~a_{i}%
\in k\}$ is a $k$-algebra isomorphic to $k[G]$. \ 
\begin{rem} The matrices in this $k$-algebra
are exactly what are called $G$-circulant matrices with entries in $k$ (with
respect to the ordering $g_{1},\cdots ,g_{n}$), alternately characterized by
the condition $A_{i,j}=A_{k,l}$ whenever $g_{i}^{-1}g_{j}=g_{k}^{-1}g_{l}$.
\ (When $G=\left\langle g\right\rangle $ is the cyclic group of order $n$
with ordering $1,g,g^{2},\cdots ,g^{n-1}$, a $G$-circulant matrix is just
what is usually called a circulant matrix.)\ See\ for example Section 2 of
\cite{[CM2]} (in which the matrix $A^g$ is denoted $P'_g$).
\end{rem}
We now prove the Maschke theorem:

\begin{proof}
For convenience, order $G$ so that $g_{1}=1$ (the identity element of $G$). \ Suppose $0\neq
\sum\limits_{i=1}^{n}a_{i}g_{i}\in rad(k[G])$. \ We can assume without loss
of generality that $a_{1}\neq 0.$ (If, say, $a_{j}\neq 0$ we can multiply by 
$g_{j}^{-1}$.) \ Let $B=\sum\limits_{i=1}^{n}a_{i}A^{g_{i}}$ be the
corresponding $G$-circulant matrix, and note that all entries of the
principal diagonal of $B$ are $a_{1}$. \ Of course, since $%
\sum\limits_{i=1}^{n}a_{i}g_{i}$ is nilpotent, so is $B$. \ But  trace$(B)=na_{1}\neq 0$ since $a_{1}\neq
0\neq n\in k$.  
 It follows that the matrix $B$ cannot be nilpotent, a
contradiction.  Thus there are no nonzero elements of $rad(k[G])$.
\end{proof}

\end{document}